\magnification=1100
\baselineskip=12pt
\hsize12cm
\vsize18cm
\font\twelverm=cmr12
\font\twelvei=cmmi12
\font\twelvesy=cmsy10
\font\twelvebf=cmbx12
\font\twelvett=cmtt12
\font\twelveit=cmti12
\font\twelvesl=cmsl12

\font\ninerm=cmr9
\font\ninei=cmmi9
\font\ninesy=cmsy9
\font\ninebf=cmbx9
\font\ninett=cmtt9
\font\nineit=cmti9
\font\ninesl=cmsl9

\font\eightrm=cmr8
\font\eighti=cmmi8
\font\eightsy=cmsy8
\font\eightbf=cmbx8
\font\eighttt=cmtt8
\font\eightit=cmti8
\font\eightsl=cmsl8

\font\sixrm=cmr6
\font\sixi=cmmi6
\font\sixsy=cmsy6
\font\sixbf=cmbx6

\catcode`@=11 
\newskip\ttglue
\def\twelvepoint{\def\rm{\fam0\twelverm}
\textfont0=\twelverm  \scriptfont0=\ninerm  
\scriptscriptfont0=\sevenrm
\textfont1=\twelvei  \scriptfont1=\ninei  \scriptscriptfont1=\seveni
\textfont2=\twelvesy  \scriptfont2=\ninesy  
\scriptscriptfont2=\sevensy
\textfont3=\tenex  \scriptfont3=\tenex  \scriptscriptfont3=\tenex
\textfont\itfam=\twelveit  \def\it{\fam\itfam\twelveit}%
\textfont\slfam=\twelvesl  \def\sl{\fam\slfam\twelvesl}%
\textfont\ttfam=\twelvett  \def\tt{\fam\ttfam\twelvett}%
\textfont\bffam=\twelvebf  \scriptfont\bffam=\ninebf
\scriptscriptfont\bffam=\sevenbf  \def\bf{\fam\bffam\twelvebf}%
\tt  \ttglue=.5em plus.25em minus.15em
\normalbaselineskip=15pt
\setbox\strutbox=\hbox{\vrule height10pt depth5pt width0pt}%
\let\sc=\tenrm  \let\big=\twelvebig  \normalbaselines\rm}

\def\tenpoint{\def\rm{\fam0\tenrm}
\textfont0=\tenrm  \scriptfont0=\sevenrm  \scriptscriptfont0=\fiverm
\textfont1=\teni  \scriptfont1=\seveni  \scriptscriptfont1=\fivei
\textfont2=\tensy  \scriptfont2=\sevensy  \scriptscriptfont2=\fivesy
\textfont3=\tenex  \scriptfont3=\tenex  \scriptscriptfont3=\tenex
\textfont\itfam=\tenit  \def\it{\fam\itfam\tenit}%
\textfont\slfam=\tensl  \def\sl{\fam\slfam\tensl}%
\textfont\ttfam=\tentt  \def\tt{\fam\ttfam\tentt}%
\textfont\bffam=\tenbf  \scriptfont\bffam=\sevenbf
\scriptscriptfont\bffam=\fivebf  \def\bf{\fam\bffam\tenbf}%
\tt  \ttglue=.5em plus.25em minus.15em
\normalbaselineskip=12pt
\setbox\strutbox=\hbox{\vrule height8.5pt depth3.5pt width0pt}%
\let\sc=\eightrm  \let\big=\tenbig  \normalbaselines\rm}

\def\ninepoint{\def\rm{\fam0\ninerm}
\textfont0=\ninerm  \scriptfont0=\sixrm  \scriptscriptfont0=\fiverm
\textfont1=\ninei  \scriptfont1=\sixi  \scriptscriptfont1=\fivei
\textfont2=\ninesy  \scriptfont2=\sixsy  \scriptscriptfont2=\fivesy
\textfont3=\tenex  \scriptfont3=\tenex  \scriptscriptfont3=\tenex
\textfont\itfam=\nineit  \def\it{\fam\itfam\nineit}%
\textfont\slfam=\ninesl  \def\sl{\fam\slfam\ninesl}%
\textfont\ttfam=\ninett  \def\tt{\fam\ttfam\ninett}%
\textfont\bffam=\ninebf  \scriptfont\bffam=\sixbf
\scriptscriptfont\bffam=\fivebf  \def\bf{\fam\bffam\ninebf}%
\tt  \ttglue=.5em plus.25em minus.15em
\normalbaselineskip=11pt
\setbox\strutbox=\hbox{\vrule height8pt depth3pt width0pt}%
\let\sc=\sevenrm  \let\big=\ninebig  \normalbaselines\rm}

\def\eightpoint{\def\rm{\fam0\eightrm}
\textfont0=\eightrm  \scriptfont0=\sixrm  \scriptscriptfont0=\fiverm
\textfont1=\eighti  \scriptfont1=\sixi  \scriptscriptfont1=\fivei
\textfont2=\eightsy  \scriptfont2=\sixsy  \scriptscriptfont2=\fivesy
\textfont3=\tenex  \scriptfont3=\tenex  \scriptscriptfont3=\tenex
\textfont\itfam=\eightit  \def\it{\fam\itfam\eightit}%
\textfont\slfam=\eightsl  \def\sl{\fam\slfam\eightsl}%
\textfont\ttfam=\eighttt  \def\tt{\fam\ttfam\eighttt}%
\textfont\bffam=\eightbf  \scriptfont\bffam=\sixbf
\scriptscriptfont\bffam=\fivebf  \def\bf{\fam\bffam\eightbf}%
\tt  \ttglue=.5em plus.25em minus.15em
\normalbaselineskip=9pt
\setbox\strutbox=\hbox{\vrule height7pt depth2pt width0pt}%
\let\sc=\sixrm  \let\big=\eightbig  \normalbaselines\rm}

\def\twelvebig#1{{\hbox{$\textfont0=\twelverm\textfont2=\twelvesy
	\left#1\vbox to10pt{}\right.\n@space$}}}
\def\tenbig#1{{\hbox{$\left#1\vbox to8.5pt{}\right.\n@space$}}}
\def\ninebig#1{{\hbox{$\textfont0=\tenrm\textfont2=\tensy
	\left#1\vbox to7.25pt{}\right.\n@space$}}}
\def\eightbig#1{{\hbox{$\textfont0=\ninerm\textfont2=\ninesy
	\left#1\vbox to6.5pt{}\right.\n@space$}}}
 
\def\tit{\bigskip}

\def\Pt{{\bf P}^3}
\def\Ptr{{\bf P}^3}
\def\Pq{{\bf P}^4}
\def\Pcq{{\bf P}^5}

\def\Pn{{\bf P}^n}
\def\epf{{\vrule height.9ex width.8ex depth-.1ex}}

\def\o{{\cal O}}
\def\i{{\cal I}}

\topinsert 
\endinsert
\font\big=cmbx10 scaled \magstep2

\centerline{\big \hbox {On smooth surfaces in projective four-space lying on}}
\par
\centerline{\big \hbox {quartic hypersurfaces with isolated singularities.}}
\vskip 0.8truecm
\centerline {Ph. Ellia $^*$ - D. Franco $^{**}$ \footnote {$^1$}
{Partially supported by MURST and Ferrara Univ. in the framework
of the project:
"Geometria algebrica, algebra commutativa e aspetti computazionali"}}
\medskip
{\eightpoint 
\centerline {Dipartimento di Matematica, Universit\`a di Ferrara}
\centerline {via Machiavelli 35 - 44100 Ferrara, Italy}
\centerline {$^{*}$ e-mail: phe@dns.unife.it}
\centerline {$^{**}$ e-mail: frv@dns.unife.it}}
\tit
\bigskip

\centerline {\it Dedicated to Robin Hartshorne in occasion of his 60th birthday.}
\bigskip

\centerline {\bf Introduction }
\tit
In the classification of smooth codimension two subvarieties of $\Pn$, surfaces in $\Pq$ (resp. threefolds in $\Pcq$) seem to lie between two extremal situations: every curve can be embedded in $\Ptr$ while, according to Hartshorne's conjecture, every smooth, codimension two subvariety of $\Pn$, $n \geq 6$, should be a complete intersection.
\par
As well known, not every surface can be embedded in $\Pq$ and, for example, the degrees of smooth rational surfaces in $\Pq$ are bounded (as conjectured by Hartshorne and Lichtenbaum some years ago); more precisely if $S \subset \Pq$ is a smooth surface not of general type, then $deg(S) \leq 46$ ([C]). By the way, it is believed that this result is not optimal and it is conjectured that the sharp bound should be $deg(S) \leq 15$.
\par
It follows from general facts that surfaces of non general type usually lie on hypersurfaces of low degree. So it seems natural to approach this problem by classifying surfaces on hypersurfaces of low degree (a question of independant interest). For hypercubics this has been done by Koelblen ([K]). In this paper we consider surfaces on hyperquartics with isolated singularities, and we prove:
\tit
{\bf Theorem.} {\it
Let $S \subset \Pq$ be a smooth surface of degree $d$ lying on a quartic hypersurface with isolated singularities.
\par \noindent
(i) If $p_g=0$ (so in particular if $S$ is rational), then $d \leq 23$.
\par \noindent
(ii) If $h^0(\omega_S(-1))=0$ (so in particular if $S$ is of non general type), then $d \leq 27$.}
\tit
In a few words, our proof (inspired by [A]) goes as follows. On one hand it is not hard, under the assumption $h^0(\omega _S(-1))=0$ (resp. $p_g=0$), to get a lower bound $h^2(\i_S(k)) \geq f(k)$ (see Cor. 7, here $d=4k+r,0 \leq r \leq 3$). On the other hand, we also have an upper bound: $h^2(\i_S(k)) \leq g(k)$ (see Cor. 3, Remark 3.1), but $g(k)$ depends on cohomological invariants of $C$, the general hyperplane section of $S$. We would like to derive a contradiction from the inequality $f(k) \leq g(k)$. In order to apply this naive plan we must have precise informations on the cohomological invariants of $C$. For this we notice that, by the "Jacobi's formula", $\pi = g(C)$ is not too far from $G(d,4)$, the maximal genus of a curve of degree $d$ on a quartic surface (see Lemma 8). Here the assumption that the hyperquartic has only isolated singularities is crucial. Then we proceed to a detailed study of the cohomological invariants of curves with genus close to $G(d,4)$ (Lemma 9, Propositions 12, 14, 16, 18).
\par
It seems worthwhile to stress that the assumptions we really need involve the degrees of the minimal generators of $H^0_*(\omega _S)$, so that our theorem apply also to some surfaces of general type.  
\tit
\centerline{\bf Generalities }
\tit

{\bf Lemma 1.}{ \it 
Let $S\subset \Pq$ be a smooth surface and let $C$ denote its general hyperplane section. 
\par
Set $c:= max \{ k/ h^1(\i_C(k))\not= 0 \}$, $v:= max \{ k/ h^2(\i_S(k)) \not= 0 \}$.
\par
We have $h^2(\i_S(t)) \leq \sum\limits_{m \geq t+1}^{v+1}h^1(\i_C(m))$ (in particular $h^2(\i_S(m))= 0$ if $m \geq c$).}

\tit
{\bf Proof.} From the exact sequence
\par \noindent
$$0 \to \i_{S}(t) \to \i_{S}(t+1) \to \i_C(t+1) \to 0$$
it follows that: $h^2(\i_S(t)) \leq h^1(\i_C(t+1)) + h^2(\i_S(t+1))$, we argue by descending induction starting from $t=v$. \epf
\tit
{\bf The numerical character of a set of points in the plane.}
\tit
To every zero-dimensional subscheme $\Gamma \subset {\bf P}^2$ there is associated a sequence of integers $\chi (\Gamma)=(n_0,...,n_{\sigma -1})$, called the numerical character of $\Gamma$, which encodes the Hilbert function of $\Gamma$. We recall the basic properties of the numerical character:
\par 
(i) $\sigma = min\{ k/ h^0(\i_{\Gamma}(k))\neq 0)\}$
\par
(ii) $\sum\limits_{i=o}^{\sigma -1}(n_i -i)=d$
\par
(iii) $h^1(\i_{\Gamma}(n))=h_{\chi}(n):=\sum\limits_{i=0}^{\sigma -1}[(n_i-n-1)_+ - (i-n-1)_+]$ (where $(x)_+ =max\{x,0\}$).
\par
The genus of a numerical character is: $g(\chi):=\sum\limits_{m \geq 1} h_{\chi}(m)$.
\par
If $C \subset \Pt$ is an integral curve, its numerical character, $\chi (C)$, is the numerical character of its general plane section.
\par
The numerical character of an integral curve is connected: $n_i - n_{i+1} \leq 1, 0 \leq i \leq \sigma -2$ (cf [GP]). Moreover $p_a(C) \leq g(\chi(C))$, with equality if and only if $C$ is a.C.M. (arithmetically Cohen-Macaulay).
\tit
Let $\rho _m: H^0(\i_C(m)) \to H^0(\i_{\Gamma}(m))$ denote the natural map of restriction ($\Gamma = C \cap H$ is the general plane section of $C$). Set $R_m:= coker(\rho _m)$ so that we have an exact sequence of modules of finite length:
\par \noindent
$$0 \to R \to M(-1) \to M \to Q \to 0$$
where $M = H^1_*(\i_C)$ and where $Q_m = ker(H^1(\i_{\Gamma}(m)) \to H^2(\i_C(m-1)))$. Finally set $r_m = dim(R_m)$ and define $q_m$ similarly.
\tit
{\bf Lemma 2.}{ \it
With notations as above:
\par \noindent
(i) $\sum\limits_{m \geq 1}r_m = \sum\limits_{m \geq 1}q_m = g(\chi(C)) - g(C)$.
\par \noindent
(ii) If ${\bf I}(\Gamma)$ is generated in degrees $\leq t_0$ and if $\rho_m$ is surjective for some $m \geq t_0$, then $\rho_t$ is surjective for $t \geq m$.
\par \noindent
(iii) $\forall m$: $h^1(\i_C(m)) \leq \sum\limits_{t \geq m+1}r_t \leq g(\chi(C)) - g(C)$.}
\tit
{\bf Proof.} (i) Follows from $g(\chi(C)) = \sum\limits_{m \geq 1}h^1(\i_{\Gamma}(m))$ and the exact sequences, for $m \geq 1$:
\par \noindent
$$ 0 \to Q_m \to H^1(\i_{\Gamma}(m)) \to H^2(\i_C(m-1)) \to H^2(\i_C(m)) \to 0$$
(ii) Clear.
\par \noindent
(iii) By descending induction. Consider the exact sequence:
\par \noindent
$$0 \to R \to M(-1) \to M \to Q \to 0$$
at level $c+1$ it yields $r_{c+1} = h^1(\i_C(c))$. Then at level $m$ it gives: 
\par \noindent
$h^1(\i_C(m-1)) \leq r_m + h^1(\i_C(m))$. \epf
\tit
{\bf Corollary 3.}{ \it
Let $S \subset \Pq$ be a smooth surface and let $C$ denote its general hyperplane section. Assume $c > e-1$, then:
\par 
$h^2(\i_S(t)) \leq [(c - t).(g(\chi(C)) - g(C))]_+$.}
\tit
{\bf Proof.} By Lemma 1: $h^2(\i_S(t)) \leq \sum\limits_{m = t+1}^c h^1(\i_C(m))$. Using Lemma 2(iii): 
\par \noindent
$$h^2(\i_S(t)) \leq \sum\limits_{m \geq t+2}r_m + \sum\limits_{m \geq t+3}r_m+ ...+ \sum\limits_{m \geq c}r_m + r_{c+1}$$
These are $c - t$ terms and each one is bounded by $\sum\limits_{m \geq 1}r_m = g(\chi(C)) - g(C)$ (see Lemma 2(i)). \epf
\tit
{\bf Remark 3.1.}{ \it
The bound of the corollary is very rough, we can improve it, for example, as follows:}
\par \noindent
$$h^2(\i_S(t)) \leq [(c - t).(g(\chi(C)) - g(C))-(c-t).\sum\limits_{m \leq t+1}r_m]_+$$
\tit
{\bf Lemma 4.}{ \it
Let $C$ be a smooth connected curve of degree $d$ in $\Ptr$ lying on a smooth surface of degree $s$. Then: $c \leq d + e(1-s) + s^2 - 4s$. }
\tit
{\bf Proof.} See [E2], lemme VI.3. \epf
\tit
{\bf Lemma 5.}{ \it
Let $S \subset \Pq$ be a smooth surface with $h^0(\omega _S(-1)) = 0$, then:
\par \noindent
(i) $p_g \leq \pi - {{d} \over {2}}$
\par \noindent
(ii) If $C$ is linearly normal (for instance if $q = 0$ and $d > 4$) then: $p_g \leq \pi - d + 3$.}
\tit
{\bf Proof.} From the assumption $h^0(\omega_S(-1))=0$ and the exact sequence:
\par \noindent
$$0 \to \omega _S(-1) \to \omega _S \to \omega _C(-1) \to 0$$
we get: $p_g \leq h^0(\omega_C(-1))$.
\par \noindent
(i) By Clifford's theorem (of course $\omega_C(-1)$ is special), $h^0(\omega_C(-1)) \leq \pi - {{d} \over {2}}$.
\par \noindent
(ii) Use $h^0(\omega_C(-1)) = h^1(\o_C(1))$. \epf
\tit
\centerline {\bf Surfaces on hyperquartics with isolated singularities.}
\tit
{\it Notation:} We denote by $G(d,s)$ the maximal genus of a curve of degree $d$ in $\Ptr$ not lying in a surface of degree $<s$. If $d > s(s-1)$ then $G(d,s)=1+{{d(d+s^2 -4s)}\over{2s}}-{{r(s-1)(s-r)}\over{2s}}$, where $d+r =0$ (mod s), $0 \leq r <s$.
\tit
{\bf Lemma 6.} {\it
Let $S \subset \Pq$ be a smooth surface of degree $d > 16$ lying on an irreducible hypersurface of degree four. Set $d = 4k + r, 0 \leq r \leq 3$ and $\pi = G(d,4) - \delta$. Then:
\par \noindent
$h^2(\i_S(k)) \geq {{2}\over {3}}k^3 + k^2({{r}\over{2}} - 1)+k({{7}\over {3}}+{{r^2}\over {2}}-2r-\delta)-p_g.$}
\tit
{\bf Proof.} We have $\chi (\i_S(k)) = \chi (\o_{\Pq}(k)) - \chi(\o_S(k))$, so $h^2(\i_S(k))=h^0(\o_{\Pq}(k)) - h^0(\i_S(k))-\chi(\o_S(k))+h^1(\i_S(k))+h^3(\i_S(k))$, hence: $h^2(\i_S(k))\geq h^0(\o_{\Pq}(k)) - h^0(\i_S(k))-\chi(\o_S(k))$. Since $S$ lies on an irreducible quartic hypersurface, and since $d > 16$, $h^0(\i_S(k))=h^0(\o_{\Pq}(k-4))$. It follows that $h^2(\i_S(k)) \geq {{2}\over{3}}k^3+k^2+{{7}\over {3}}k+1-\chi(\o_S(k))$. By Riemann-Roch: $\chi(\o_S(k)) = {{kH.(kH-K)} \over {2}}+\chi = {{d(k+1)k}\over {2}}-k(\pi - 1)+1-q+p_g$. So $h^2(\i_S(k)) \geq {{2}\over{3}}k^3+k^2+{{7}\over {3}}k-{{d(k+1)k} \over {2}}+k(\pi - 1) - p_g$. Taking into account that: $d = 4k+r$ and that $\pi = G(d,4) - \delta = 1+2k^2+kr+{{r} \over {2}}(r-3)-\delta$, we get the result after a little computation. \epf
\tit
{\bf Corollary 7.} {\it
With notations as above:
\par \noindent
(i) If $p_g = 0$ then $h^2(\i_S(k)) \geq \rho_{\delta,r}(k)$ where:
\par
$\rho_{\delta,r}(k) = {{2} \over {3}}k^3+k^2({{r} \over {2}}-1)+k({{7}\over{3}}+{{r^2}\over{2}}-2r-\delta)$.
\par \noindent
(ii) If $h^0(\omega _S(-1))=0$ and if the general hyperplane section of $S$ is linearly normal in $\Ptr$, then: $h^2(\i_S(k)) \geq \lambda_{\delta,r}(k)$, where:
\par
$\lambda_{\delta,r}(k)={{2} \over {3}}k^3+k^2({{r} \over {2}}-3)+k({{19}\over {3}}+{{r^2}\over {2}}-3r-\delta)-{{r(r-5)}\over {2}}-4+\delta$.
\par \noindent
(iii) If $h^0(\omega_S(-1))=0$ then: $h^2(\i_S(k)) \geq \phi_{\delta,r}(k)$, where:
\par
$\phi_{\delta,r}(k)={{2} \over {3}}k^3+k^2({{r} \over {2}}-3)+k({{13}\over {3}}+{{r^2}\over {2}}-3r-\delta)+2r+\delta-1-{{r^2}\over {2}}$.}
\tit
{\bf Proof.} (i) Follows directly from Lemma 6.
\par \noindent
(ii) We use $p_g \leq \pi -d+3$, see Lemma 5.
\par \noindent
(iii) We use $p_g \leq \pi - {{d}\over {2}}$. \epf
\tit
{\bf Remark 7.1.} {\it
For later use we observe that $\phi_{\delta, r}(k), \lambda_{\delta,r}(k), \rho_{\delta,r}(k)$ are increasing functions of $k$ for $0 \leq \delta \leq 10$ and $0 \leq r \leq 3$.}
\tit
{\bf Lemma 8.} {\it
Let $S \subset \Pq$ be a smooth surface of degree $d$ lying on a hypersurface of degree four with isolated singularities. Then $\pi = G(d,4)-\delta$ and, if $d=4k+r, 0 \leq r \leq 3$, then: $\delta \leq 10$ if $r=0$, $\delta \leq 9$ if $r=1$ or $r=3$; $\delta \leq 8$ if $r=2$.}
\tit
{\bf Proof.} This follows from the "Jacobi's formula": $\pi=1+{{d^2-\mu}\over {8}}$ (see [A], lemma 2.1). Since the hypersurface has only isolated singularities: $\mu \leq 81$, and since $G(d,4)=1+{{d^2 -3r(4-r)}\over {8}}$, we get $\delta = {{\mu -3r(4-r)}\over {8}}$ with $\mu \leq 81$ and we conclude. \epf
\tit
{\it Arithmetically Cohen-Macaulay surfaces.}
\tit
{\bf Lemma 9.} {\it
Let $S \subset \Pq$ be a smooth surface of degree $d$ lying on an irreducible hypersurface of degree four. If $S$ is arithmetically Cohen-Macaulay and if $h^0(\omega _S(-1))=0$ (resp. $p_g=0$), then $d \leq 16$ (resp. $d \leq 12$).}
\tit
{\bf Proof.} Since $S$ is a.C.M., the minimal free resolution yields an exact sequence:
$$
0 \to \bigoplus\limits^r_{j=1} \o(-b_j) \to \bigoplus\limits^{r+1}_{i=1} \o(-a_i) \to \i_S \to 0
$$
Dualizing, we get:
$$
0 \to \o(-6) \to \bigoplus\limits^{r+1}_{i=1} \o(a_i-6) \to \bigoplus\limits^r_{j=1} \o(b_j-6) \to \omega_S(-1) \to 0
$$
If $h^0(\omega_S(-1))=0$ then $b^+ \leq 5$ (here $b^+=max\{b_j\}$). Since $b^+>a^+$, we get $a_i \leq 4$, $\forall i$. Since $r \geq 2$, $h^0(\i_S(4)) \geq 2$ and this implies, by Bezout, $d \leq 16$.
\par
If $h^0(\omega_S)=0$, the proof is similar. \epf
\tit
{\bf Remark 9.1.} Thanks to the previous lemma, in the sequel, we will assume $S$ not a.C.M., i.e. we will assume that the general hyperplane section, $C \subset \Pt$ of $S$ is not projectively normal.
\tit
\centerline {\bf The case $d=4k$.}
\tit
{\bf Lemma 10.} {\it
Let $C \subset \Ptr$ be a smooth, connected curve of degree $d = 4k$, $k \geq 4$ and genus $\pi$, lying on a smooth quartic surface. Assume $\pi = G(d,4) - \delta$ with $\delta \leq 10$. Moreover suppose $C$ non projectively normal. Then:
\par \noindent
(i) $g(\chi (C)) = G(d,4) - 2.$
\par \noindent
(ii) $\delta \geq 3$ and $c > e-1.$
\par \noindent
(iii) $k-3 \leq e \leq k-1$. Moreover if $e = k-3$ then $\delta = 10$, $h^1(\i_C(k-2))=0, C$ is linearly normal and: $c \leq k+3$ or $r_{k+1} = 2$.}
\tit
{\bf Proof.} (i) There are only two connected numerical characters of degree $4k$, length $4$: $\phi = (k+3, k+2, k+1, k)$, $\chi = (k+2, k+2, k+1, k+1)$. If $\chi (C)=\phi$ then by [D] (see also [E]), then $C$ is projectively normal, which is excluded. Hence $\chi(C)=\chi$. We have: $h_{\phi}(k)=3, h_{\phi}(k+1)=1$, $h_{\chi}(k)=2, h_{\chi}(k+1)=0$ while $h_{\phi}(t)=h_{\chi}(t)$ otherwise. It follows that $g(\chi)=g(\phi)-2=G(d,4)-2$. 
\par \noindent
(ii) Since $C$ is not projectively normal, we have $g(C) < g(\chi(C))$, i.e. $\delta \geq 3$. Since $dk = 4k^2 > 2\pi -2=4k^2-2\delta$, we get $e \leq k-1$. Since $C$ lies on an irreducible quartic surface and since $d > 16$, $h^0(\i_C(t))=h^0(\o_{\Ptr}(t-4))$ for $t \leq k$, and a simple computation shows that: $h^1(\i_C(k))=\delta - 2 \neq 0$, hence $c \geq k$. It follows that $c > e-1$.
\par \noindent
(iii) We have seen (cf (ii)) that $e \leq k-1$. Since $C$ lies on an irreducible quartic surface, for $t \leq k, h^1(\i_C(t))=dt-\pi +1 + h^1(\o_C(t)) + h^0(\o_{\Ptr}(t-4))-h^0(\o_{\Ptr}(t))$. So if $e < k-2$, $h^1(\i_C(k-2))=\delta -10$, hence the only possibility is $\delta = 10$ and $h^1(\i_C(k-2))=0$, this implies also $e = k-3$ (otherwise $C$ would be projectively normal by Castelnuovo-Mumford's lemma). In this case, since $h^1(\i_C(k-2))=0$, by descending induction, $h^1(\i_C(t))=0$ if $t \leq k-2$, in particular $C$ is linearly normal.
\par \noindent
If ${\bf I}(C)$ has a generator of degree $k+1$, we can link $C$ to a curve, $X$, by a complete intersection $(4,k+1)$. The curve $X$ has degree $4$ and $p_a(X)=-7$: since $X$ lies on a smooth quartic surface the only possibility is the disjoint union of two double lines (of arithmetic genus $-3$), see lemma below. In particular $h^1(\i_X(-t))=0$ if $t > 2$. By liaison this implies $c \leq k+3$.
\par \noindent
If ${\bf I}(C)$ has no generator of degree $k+1$, then $r_{k+1}= 2$. \epf
\tit
{\bf Lemma 11.} {\it
Let $X$ be a curve of degree $4$ and arithmetic genus $-7$, lying on a smooth quartic surface. Then $X$ is the disjoint union of two double lines of arithmetic genus $-3$.}
\tit
{\bf Proof.} We have $X^2 = 2p-2 = -16$. Clearly $X$ is non-reduced, so it must contain a line with multiplicity $\geq 2$ or $X$ is a double conic. Checking case by case we get the lemma. \epf
\tit
{\bf Proposition 12.} {\it
Let $S \subset \Pq$ be a smooth surface of degree $d = 4k, k \geq 5$ lying on a hyperquartic with isolated singularities.
\par \noindent
(i) If $h^0(\omega _S(-1))=0$, then $d \leq 24$.
\par \noindent
(ii) If $p_g =0$, then $d \leq 20$.}
\tit
{\bf Proof.} We have $\pi = G(d,4)-\delta$ and from Lemma 8, we may assume $\delta \leq 10$. From Lemma 10 (see Remark 9.1): $k-3 \leq e \leq k-1$ and $c > e-1$. First we make the following:
\tit
{\it Claim:} $h^2(\i_S(k)) \leq 6\delta - 12$ if $e \geq k-2$ and, if $e = k-3$ then $C$ is linearly normal, $\delta = 10$, and $h^2(\i_S(k)) \leq 54$.
\tit
Assume this for a while and let's conclude the proof.
\par \noindent
(i) If $e=k-3$, since $C$ is linearly normal, we have (see Corollary 7 and Remark 7.1): $h^2(\i_S(k)) \geq \lambda_{\delta,0}(7)=122-6\delta$, if $k \geq 7$. From the claim it follows that: $54 \geq 122 - 6\delta$, which is impossible if $\delta = 10$.
\par \noindent
Now assume $e \geq k-2$. This time we use $h^2(\i_S(k)) \geq \phi_{\delta,0}(k)$. Arguing as above: $h^2(\i_S(k)) \geq \phi_{\delta,0}(7) = 111 - 6\delta$. Combining with the claim: $123 \leq 12\delta$, which is impossible for $\delta \leq 10$. It follows that it must be $k \leq 6$, i.e. $d \leq 24$.
\par \noindent
(ii) We argue as before but using $h^2(\i_S(k)) \geq \rho_{\delta,0}(k) \geq \rho_{\delta,0}(6) = 122 - 6\delta$. If $e \geq k-2$, we get $122 - 6\delta \leq 6\delta -12$ which is impossible for $\delta \leq 10$. If $e = k-3$, since $\delta = 10$, we have $h^2(\i_S(k)) \geq 62$ and we conclude with the claim.
\tit
To conclude, let's prove the claim:
\par \noindent
We have $h^2(\i_S(k)) \leq (c-k).(g(\chi(C)) - \pi - \sum\limits_{m \leq k+1} r_m)$, $(*)$, (see Remark 3.1). Since $g(\chi(C))=G(d,4)-2$, it follows that: $h^2(\i_S(k)) \leq (c-k).(\delta - 2)$, $(**)$. By Lemma 4: $c \leq d - 3e = 4k - 3e$, $(***)$. 
\par \noindent
If $e \geq k-2$, from $(**)$ and $(***)$ we get: $h^2(\i_S(k)) \leq 6\delta - 12$ and the claim is proved.
\par \noindent
If $e = k-3$, by Lemma 10, $\delta = 10,h^1(\i_C(t))=0$ if $t \leq k-2$, and $c \leq k+3$ or $r_{k+1}=2$.
\par \noindent
- If $c \leq k+3$, we get $h^2(\i_S(k)) \leq 3\delta - 6 = 24$.
\par \noindent
- In any case, from Lemma 4, $c \leq k+9$. Since $r_{k+1}=2$, from $(*)$, we get: $h^2(\i_S(k)) \leq 9(\delta - 4)=54$, since $\delta = 10$. \epf
\tit
\centerline {\bf The case $d = 4k+1$.}
\tit
{\bf Lemma 13.} {\it
Let $C \subset \Ptr$ be a smooth, connected curve of degree $d = 4k+1$, $k \geq 4$, and genus $\pi = G(d,4) - \delta = 2k^2+k-\delta$. Assume $C$ lies on an irreducible quartic surface and $C$ not projectively normal. Then, if $\delta \leq 9$:
\par \noindent
(i) $\delta \geq 2$ and $g(\chi(C))=G(d,4)-1$.
\par \noindent
(ii) $k-2 \leq e \leq k$.
\par \noindent
(iii) $c > e-1$.}
\tit
{\bf Proof.} (i) There are two possible connected numerical characters: the maximal one, $\phi = (k+3,k+2,k+1,k+1)$ and $\chi = (k+2,k+2,k+2,k+1)$. If $\chi(C) = \phi$ then $C$ is projectively normal ([D], [E]), so we may assume $\chi(C) = \chi$. We have $h_{\phi}(k+1)=1$, $h_{\chi}(k+1)=0$, while $h_{\phi}(t)=h_{\chi}(t)$ if $t \neq k+1$. It follows that $g(\chi)=g(\phi)-1=G(d,4)-1$. Since $C$ is not projectively normal, it must be $g(C) < g(\chi(C))$, i.e. $\delta \geq 2$.
\par \noindent
(ii) Since $d(k+1) > 2\pi -2$, we have $e \leq k$. 
\par \noindent
Assume $e < k-2$. Then $h^0(\o_C(k-2))=d(k-2)-\pi+1 =2k^2-8k-1+\delta$. On the other hand $h^0(\i_C(k-2))=h^0(\o_{\Ptr}(k-6))$ since $C$ is contained in an irreducible quartic surface. Now we must have:
\par \noindent
$h^0(\o_{\Ptr}(k-2)) - h^0(\i_C(k-2)) \leq h^0(\o_C(k-2))$, which is equivalent to: $2k^2-8k+10 \leq 2k^2-8k-1+\delta$, but this is impossible if $\delta \leq 9$.
\par \noindent
(iii) We have $h^1(\i_C(k))=h^0(\o_C(k))+h^0(\i_C(k))-h^0(\o_{\Ptr}(k))$. We have $h^0(\i_C(k))=h^0(\o_{\Ptr}(k-4))$, we get $h^1(\i_C(k))=\delta-1+h^1(\o_C(k))$, since $\delta \geq 2$ by (i), it follows that $c \geq k$ hence $c > e-1$. \epf
\tit
{\bf Proposition 14.} {\it
Let $S \subset \Pq$ be a smooth surface of degree $d = 4k+1$, with $k \geq 4$. Assume $S$ lies on an irreducible hypersurface of degree four with isolated singularities.
\par \noindent
(i) If $p_g =0$, then $d \leq 21$.
\par \noindent
(ii) If $h^0(\omega_S(-1))=0$, then $d \leq 25$.}
\tit
{\bf Proof.} (i) We have $\pi = G(d,4)- \delta$ and, by Lemma 8, we may assume $\delta \leq 9$. We may assume $C$ not projectively normal (Remark 9.1). By Lemma 13, $k-2 \leq e \leq k$ and $c >e-1$. By Corollary 7 and Remark 7.1, we have $h^2(\i_S(k)) \geq \rho_{\delta,1}(6) = 131 -6\delta$. On the other hand, by Corollary 3, $h^2(\i_S(k)) \leq (c-k).(g(\chi(C)) - \pi)$. Since, by Lemma 4, $c \leq k+7$ and since $g(\chi(C)) \leq G(d,4)-1$, by the previous lemma, it follows that $h^2(\i_S(k)) \leq 7(\delta -1)$. Combining with the previous inequality yields a contradiction if $k \geq 6$.
\par \noindent
(ii) We argue as above but using $h^2(\i_S(k)) \geq \phi_{\delta,1}(7) = 120 -{1 \over 2} - 6\delta$. We get a contradiction if $k \geq 7$. \epf
\tit
\centerline {\bf The case $d = 4k+2$.}
\tit
{\bf Lemma 15.} {\it
Let $C \subset \Ptr$ be a smooth, connected curve of degree $d = 4k+2$ with $k \geq 4$ and genus $\pi = G(d,4)-\delta =2k^2+2k-\delta$, lying on an irreducible quartic surface. Assume $C$ non projectively normal, then, if $\delta \leq 8$:
\par \noindent
(i) $k-2 \leq e \leq k$ and $g(\chi(C)) \leq G(d,4)$.
\par \noindent
(ii) $c > e-1$.
\par \noindent
(iii) $c \leq k+8$.}
\tit
{\bf Proof.} (i) Since $d(k+1) > 2\pi -2$, $e \leq k$.
\par \noindent
Since $C$ lies on an irreducible quartic surface $h^0(\i_C(t))=h^0(\o_{\Ptr}(t-4))$ for $t \leq k$. It follows that $h^1(\i_C(t))= dt-\pi+1+h^1(\o_C(t))+h^0(\o_{\Ptr}(t-4))-h^0(\o_{\Ptr}(t))$. In particular $h^1(\i_C(k-2))=\delta - 13 + h^1(\o_C(k-2))$. Since $\delta \leq 8$, this implies $e \geq k-2$.
\par \noindent
Finally it is clear that $g(\chi(C)) \leq G(d,4)$ (by the way notice that in this case [D], [E] do not apply).
\par \noindent
(ii) We have $h^1(\i_C(k)) = \delta -1+h^1(\o_C(k))$ and $h^1(\i_C(k-1)) = \delta -5+h^1(\o_C(k-1))$. If $h^1(\i_C(k)) \neq 0$, then $c > e-1$ and we are done. Assume $h^1(\i_C(k))=0$, then $\delta = 1$ and $e = k-1$. If $h^1(\i_C(k-1)) \neq 0$, $c > e-1$. Finally, if $h^1(\i_C(k-1))=0$, observe that, since for degree reasons the exact sequence:
\par \noindent
$0 \to \i_C(t-1) \to \i_C(t) \to \i_{C \cap H}(t) \to 0$ is exact on global sections for $t \leq k$, by descending induction we get $h^1(\i_C(t))=0$ for $t \leq k$. Since by assumption, $C$ is not projectively normal, this implies $c > k$ and the condition $c > e-1$ is satisfied.
\par \noindent
(iii) Follows from (i) and Lemma 4. \epf
\tit
{\bf Proposition 16.} {\it
Let $S \subset \Pq$ be a smooth surface of degree $d = 4k+2$, $k \geq 4$, lying on a quartic hypersurface with isolated singularities.
\par \noindent
(i) If $p_g=0$, then $d \leq 22$.
\par \noindent
(ii) If $h^0(\omega_S(-1))=0$, then $d \leq 26$.}
\tit
{\bf Proof.} (i) We have $\pi = G(d,4)-\delta$ and we may assume $\delta \leq 8$ (Lemma 8). We may assume $C$ not projectively normal (Remark 9.1). By the previous lemma, $k-2 \leq e \leq k$ and $c > e-1$. By Corollary 7 and Remark 7.1, we have $h^2(\i_S(k)) \geq \rho_{\delta,2}(6)= 146 - 6\delta$. By Corollary 3, $h^2(\i_S(k)) \leq (c-k).(g(\chi(C))-\pi)$. Since $c \leq k+8$ and $g(\chi(C)) \leq G(d,4)$ (Lemma 15), we get: $h^2(\i_S(k)) \leq 8\delta$. Combining everything, we get a contradiction if $k \geq 6$.
\par \noindent
(ii) We argue as above, but using $h^2(\i_S(k)) \geq \phi_{\delta,2}(7)=134-6\delta$. \epf
\tit
\centerline {\bf The case $d = 4k+3$.}
\tit
{\bf Lemma 17.} {\it
Let $C \subset \Ptr$ be a smooth, connected curve of degree $d = 4k+3, k \geq 3$ and genus $\pi = G(d,4)-\delta = 2k^2+3k+1-\delta$. Assume $C$ lies on a irreducible quartic surface and that $C$ is not projectively normal. Then, if $\delta \leq 9$:
\par \noindent
(i) $g(\chi(C)) = G(d,4)-1$ and $\delta \geq 2$.
\par \noindent
(ii) $k-2 \leq e \leq k$.
\par \noindent
(iii) $c \geq k+1$, in particular $c > e-1$.}
\tit
{\bf Proof.} (i) There are two possible connected numerical characters: the maximal one, $\phi = (k+3, k+3, k+2, k+1)$ and $\chi = (k+3, k+2, k+2, k+2)$. Since $C$ is not projectively normal, by [D], [E], $\chi(C) = \chi$. We have $h_{\phi}(k+1) = 2, h_{\chi}(k+1)=1$ and $h_{\phi}(t) = h_{\chi}(t)$ otherwise. This shows that $g(\chi)=g(\phi)-1=G(d,4)-1$. Taking into account that $g(C) < g(\chi(C))$ because $C$ is not projectively normal, we get $\delta \geq 2$.
\par \noindent
(ii) Since $d(k+1) > 2\pi-2$, $e \leq k$.
\par \noindent
If $e < k-2$ then $h^0(\o_C(k-2))=d(k-2)-\pi + 1$. Since $C$ is contained in a irreducible quartic surface, $h^0(\i_C(k-2))=h^0(\o_{\Ptr}(k-6))$, and we get $h^1(\i_C(k-2))=\delta - 16$, which is absurd if $\delta \leq 9$.
\par \noindent
(iii) We have $h^0(\i_C(k+1))=h^0(\o_{\Ptr}(k-3))$ (otherwise $C$ would be linked to a line by a complete intersection $(4, k+1)$ and thus $C$ would be projectively normal). It follows (since $e \leq k$, see (ii)) that $h^1(\i_C(k+1))=\delta - 1 \neq 0$ (because $\delta \geq 2$, see (i)), hence $c \geq k+1$ and in particular $c > e-1$. \epf
\tit
{\bf Proposition 18.} {\it
Let $S \subset \Pq$ be a smooth surface of degree $d = 4k+3$ lying on a quartic hypersurface with isolated singularities.
\par \noindent
(i) If $p_g=0$, then $d \leq 23$.
\par \noindent
(ii) If $h^0(\omega_S(-1))=0$, then $d \leq 27$.}
\tit
{\bf Proof.} (i) By Lemma 8, $\pi = G(d,4)-\delta$ with $\delta \leq 9$. By Corollary 3 we have $h^2(\i_S(k)) \leq (c-k).(g(\chi(C) - \pi)$. We may assume $C$ not projectively normal (Remark 9.1). Since $e \geq k-2$ by the previous lemma, by Lemma 4, we get $c \leq k+9$. Finally, since $g(\chi(C))=G(d,4)-1$, it follows that $h^2(\i_S(k)) \leq 9(\delta - 1)$ (*). On the other hand, by Corollary 7 and Remark 7.1, if $k \geq 6$, $h^2(\i_S(k)) \geq \rho_{\delta,3}(6)= 167-6\delta$. Combining with (*): $15\delta \geq 176$ which implies $\delta \geq 12$. So in our case, it must be $k \leq 5$ i.e. $d \leq 23$.
\par \noindent
(ii)We argue as above but using $h^2(\i_S(k)) \geq \phi_{\delta,3}(7)=155-{1\over 2}-6\delta$ if $k \geq 7$. (See Corollary 7) Combining with (*): $15\delta \geq 164-{1\over 2}$, which implies $\delta \geq 11$. It follows that under our assumptions, we must have $k \leq 6$. \epf
\tit
\centerline {\bf Conclusion.}
\tit
Gathering everything together (Lemma 9, Prop. 12, 14, 16 and 18):
\tit
{\bf Theorem 19.} {\it
Let $S \subset \Pq$ be a smooth surface of degree $d = 4k+r$, $0 \leq r \leq 3$, lying on a quartic hypersurface with isolated singularities.
\par \noindent
(i) If $p_g=0$, then $k \leq 5$ (in particular $d \leq 23$).
\par \noindent
(ii) If $h^0(\omega_S(-1))=0$, then $k \leq 6$ (in particular $d \leq 27$).}
\tit
{\bf Corollary 20.} {\it
Let $S \subset \Pq$ be a smooth surface lying on a quartic hypersurface with isolated singularities.
\par \noindent
(i) If $S$ is not of general type, then $d \leq 27$.
\par \noindent
(ii) If $S$ is rational, then $d \leq 23$.}
\tit
{\bf Proof.} Just observe that if $S$ is not of general type then $h^0(\omega_S(-1))=0$\epf
\tit
Another immediate consequence:
\tit
{\bf Corollary 21.} {\it
Let $V \subset {\bf P}^5$ be a smooth threefold of degree $d$, lying on a quartic hypersurface, $\Sigma$. Assume $h^0(\omega _V)=0=h^1(\omega _V(-1))$. \par
If $dim(Sing(\Sigma)) \leq 1$, then $d \leq 27$.}
\tit
{\bf Proof.} Let $S$ be a general hyperplane section of $V$. From the exact sequence: 
$$
0 \to \omega _V(-1) \to \omega _V \to \omega _S(-1) \to 0 $$
we get $h^0(\omega _S(-1))=0$. Then, since $S$ lies on an hyperquartic with isolated singularities, apply Theorem 19. \epf

\tit
\centerline {\bf References}
\tit
\item{[A]} Aure, A.: {\it "The smooth surfaces on cubic hypersurface in $\bf P ^4$ with isolated singularities"}. Math. Scand. {\bf 67}, No.2, 215-222 (1990)
\item{[C]} Cook, M.: {\it "A smooth surface of degree in $\bf P ^4$ not of general type has degree at most $46$"}, preprint. 
\item{[D]} Dolcetti, A.: {\it "Halphen's gaps for space curves of submaximum genus"}. Bull. Soc. Math. France, {\bf 116}, 157-170 (1988)
\item{[E]} Ellia, Ph.: {\it "Sur les lacunes d'Halphen"}, in Lect. Notes Math. {\bf 1389}, 43-65 (1989) Springer-Verlag.
\item{[E2]} Ellia, Ph.: {\it "D'autres composantes non r\'eduites de $Hilb(\Ptr)$"}. Math. Annalen, {\bf 277}, 433-446 (1987)
\item{[GP]} Gruson, L-Peskine, Ch.:{\it "Genre des courbes de l'espace projectif"}, in Lect. Notes Math., {\bf 687}, 31-60 (1978)
\item{[K]} Koelblen, L.: {\it "Surfaces de $\bf P ^4$ trac\'ees sur une hypersurface cubique"}. 
J. Reine Angew. Math. {\bf 433}, 113-141 (1992)

\end